\documentclass[smallextended]{svjour3}       
\smartqed  

\usepackage{amsmath,amssymb}
\usepackage[matrix,arrow,curve]{xy}

\newcommand{\proj}{\mathop{\mathrm{pr}}}
\newcommand{\coker}{\mathop{\mathrm{coker}}}

\journalname{Applied Categorical Structures}

\begin{document}

\title{Relative cohomology of algebraic theories}
\author{Simeon Pol'shin} 
\institute{ S. Pol'shin \at Institute for Theoretical Physics, NSC Kharkov Institute of Physics and Technology, Akademicheskaia St. 1, 61108 Kharkov, Ukraine \\
\email{polshin.s@gmail.com}}

\date{}

\maketitle

\begin{abstract}
We construct relative abelian categories in the sense of MacLane for models of algebraic systems in (co)complete abelian categories. As an example, we consider an analogue of Hochschild-Mitchell cohomology for the functor of Yoneda embedding.
\end{abstract}

\section{Introduction}

Recently, cohomology of algebraic systems was studied using the methods of algebraic topology~\cite{Schw},  Andr\'e-Quillen cohomology~\cite{Pir,Blanc}, Baues-Wirsching cohomology of small categories~\cite{Pir}, cotriple cohomology associated with a neglecting funtor~\cite{JP91} and ordinary homological algebra~\cite{JP91,Pir}. Also in~\cite{Pir,JP91} some of these approaches were compared. In the present paper we propose the relative homological algebra approach to the cohomology of multi-sorted algbraic systems.  We identify the category of models of an algebraic system in a given abelian category $\mathcal{A}$ as the category of functors from an appropriate small category to $\mathcal{A}$. Then we construct a relative abelian category in the sense of~\cite{H} starting from the same neglecting functor as in~\cite{JP91}, and we study some properties of the associated cohomology theory. Also we establish a connection between this relative cohomology and the "absolute" one, thus generalizing results of~\cite{JP91}.

After some categorical preliminaries in Sec.~2, we consider resolvents in relative pair of abelian categories in the sense of~\cite{H} for models of algebraic theories in abelian categories in Sec.~3.  As an example, in Sec.~4 we construct an analogue of Hochschild-Mitchell resolvent for small categoies~\cite{Oberst,Mitch,Cib}. This reduces the cohomology of algebraic theories to the one of locally finite (dimensional) categories considered by many authors~\cite{Xu,Hersc,Chites}. For other approaches to cohomology of algebraic theories see~\cite{Schw,Pir,Blanc}.

\section{Categorical preliminaries}

We fix some notations from category theory (see~\cite{CW,Barr} for details). We denote objects of $\mathcal{C}$ by $c,{c'},d,\ldots$, morphisms of $\mathcal{C}$ by $f,g,\ldots$, functors from $\mathcal{C}$ to $\mathcal{A}$ by $K,T,S,\ldots$ and natural transformations between these functors by $\varphi, \varkappa,\ldots$. The last term means that for any object $c$ of $\mathcal{C}$ there exists a morphism $\varphi_c$ in $\mathcal{A}$  such that the following diagram
\begin{equation}\label{phi-c-d}
	\xymatrix{
T(c) \ar[r]^{\varphi_c} \ar[d]_{T(f)} & S(c) \ar[d]^{S(f)} \\
T(d) \ar[r]_{\varphi_d} & S(d)
}
\end{equation}
commutes for any morphism $f:\ c\rightarrow d$ in $\mathcal{C}$. We write $\varphi:\ T\rightharpoonup S$ for the above natural transformation, and morphisms $\varphi_c$ are called components of $\varphi$.

We write 
 $\xymatrix{
	 \bullet \ar[r]_f\ar@/^/[rr]^{fg} & \bullet \ar[r]_g & \bullet 
}$
for composition of two morphisms in $\mathcal{C}$. 

We say that ``a'' diagram of product of two objects exists in $\mathcal{C}$ if (a) the following diagram
$$\xymatrix{
c_1 & c_1 \times c_2 \ar[l]_-{\proj_1} \ar[r]^-{\proj_2} & c_2
 }$$
exists in $\mathcal{C}$ and (b) this diagram is universal. We say that the category has binary products if the above diagram  exists for any pair of objects. It is well known that $\times$ is a bifunctor, so for any pair $(f_1 ,f_2)$ of morphisms in $\mathcal{C}$ there exist ``a'' morphism $f_1 \times f_2$ in $\mathcal{C}$  such that the diagram
\begin{equation} \label{f-times}
	\xymatrix{ c_1 \ar[r]^-{f_1} & d_1 \\
	c_1 \times c_2 \ar[r]^-{f_1 \times f_2} \ar[u]^-{\proj_1} \ar[d]^-{\proj_2} & d_1 \times d_2 \ar[u]^-{\proj_1} \ar[d]^-{\proj_1} \\
	c_2 \ar[r]^-{f_2} & d_2
	}
\end{equation}
commutes. 

Suppose a category has binary product and ``nullary'' product (i.e. terminal object), then we say that it has finite products. Suppose both $\mathcal{C}$ and $\mathcal{A}$ have finite products, then we say that a functor $F\in \mathop{\mathrm{Funct}} (\mathcal{C},\mathcal{A})$ is product preserving if it takes diagrams of finite products in $\mathcal{C}$ to those in $\mathcal{A}$. It is easily seen that product preserving functors and natural transformations between them form a full subcategory $\mathop{\mathrm{FPFunct}} (\mathcal{C},\mathcal{A})$ in $\mathop{\mathrm{Funct}} (\mathcal{C},\mathcal{A})$. Let $\varphi:\ T\rightharpoonup S$ be a natural transformation between two product-preserving functors, then comparing~\eqref{phi-c-d} and~\eqref{f-times} we obtain
\begin{equation}\label{phi-times}
	\varphi_{c_1 \times c_2} \simeq \varphi_{c_1} \times \varphi_{c_2}.
\end{equation}
 
We call a category preadditive if it is enriched in $\mathbf{Ab}$ (some authors use the term ``additive'' in this case). We say that a functor is preadditive if it is an $\mathbf{Ab}$-functor in the sense of enriched category theory, and we retain the notations  $\mathop{\mathrm{Funct}} (\mathcal{C},\mathcal{A})$ and $\mathop{\mathrm{FPFunct}} (\mathcal{C},\mathcal{A})$ for categories of preadditive and preadditive product-preserving functors respectively. If a preadditive category has a terminal object, then it is also a zero one, and we denote the category of normalized (i.e. zero object-preserving) functors by $\mathop{\mathrm{Funct}_\mathrm{N}} (\mathcal{C},\mathcal{A})$.

Let $\mathcal{A}$ be a preadditive category, then a diagram of the form
\begin{equation}\label{dirsum-c}
	\xymatrix{
a_1 \ar@<1ex>[r]^-{\iota_1} & a_1 \oplus a_2 \ar@<1ex>[l]^-{\pi_1} \ar@<-1ex>[r]_-{\pi_2} & a_2 \ar@<-1ex>[l]_-{\iota_2}
 }
\end{equation}
is called diagram of direct sum if the equalities
$$\pi_1 \iota_1 + \pi_2 \iota_2 =1_{c_1 \oplus c_2} , 
\qquad \iota_k \pi_l = \left\{ 
\begin{aligned}
1_{c_k}, \ k=l \\
0, \ k\not=l
\end{aligned}
 \right. $$
are satisfied. It is well known that a diagram of direct sum is universal whenever it exists, so it defines a diagram of products in $\mathcal{A}$ . Conversely, any  diagram of product in preadditive category defines the one of direct sum, and preadditive category  which has direct sums of any pair objects and has ``a'' zero object is called additive category, so an additive category always has finite products.

Let $\mathcal{C}$ be category with finite products and let $\mathcal{A}$ be an additive category. Define objectwise direct sum $F\oplus G$ of two product preserving functors $F,G$ from $\mathcal{C}$ to $\mathcal{A}$, then considering an iterated direct sum in $\mathcal{A}$ we see that $F\oplus G$ is in turn a product-preserving functor, so $\mathop{\mathrm{FPFunct}} (\mathcal{C},\mathcal{A})$ is an additive category.

	Let $\mathcal{C}$ be a category, denote by $\mathbb{Z}\mathcal{C}$ a preadditive category with $\mathop{Ob}\mathbb{Z}\mathcal{C}=\mathop{Ob} \mathcal{C}$ and $\mathop{Mor}\mathbb{Z}\mathcal{C}=\mathbb{Z}\mathop{Mor} \mathcal{C}$ where $\mathbb{Z}$ is the functor of free abelian group, then we have an identification
	\begin{equation}\label{fp-zn}
\mathop{\mathrm{Funct}} (\mathcal{C},\mathcal{A}) \cong \mathop{\mathrm{Funct}} (\mathbb{Z}\mathcal{C},\mathcal{A})  
\end{equation}
provided $\mathcal{A}$ is a preadditive category, so we can consider $\mathop{\mathrm{FPFunct}} (\mathcal{C},\mathcal{A})$ as a subcategory of $\mathop{\mathrm{Funct}} (\mathcal{C},\mathcal{A})$ provided $\mathcal{C}$ has finite products and $\mathcal{A}$ is additive,
so we can use the theory of preadditive functors developed in~\cite{Mitch}.

\section{Relative abelian categories}

An additive category $\mathcal{A}$ is called abelian if it  obeys the following three axioms:

\begin{enumerate}
	\item[({\rm Abel-1})] $\ker f$ and $\coker f$ are nonempty for any morphism $f$ of $\mathcal{A}$.
	\item[({\rm Abel-2})] Suppose $f$ is a monomorphism and $g$ is an epimorphism, then $f\in\ker g$ iff $g\in \coker f$.
	\item[({\rm Abel-3})] We can decompose any morphism in $\mathcal{A}$ into the epimorphism followed by monomorphism.
\end{enumerate}

It is well known that the decomposition provided by (Abel-3) is functorial.

Let $T,S$ be product preserving functors from $\mathcal{C}$ to $\mathcal{A}$ and let $\varphi:\ T\rightharpoonup S$  be a natural transformation. Like the case of ordinary functors~\cite{Toh,H}, for any object $c$ of $\mathcal{C}$ we can choose a morphism ${\varkappa_c}:\ K(c)\rightarrow C(c)$ in $\mathcal{A}$ such that $\varkappa_c \in\ker \varphi_c$. Like the case of ordinary functors, it follows that an exact sequence
$$\xymatrix{
0\ar[r] & K(c) \ar[r]^{\varkappa_c} & T(c) \ar[r]^{\varphi_c} & S(c)
}$$
is functorial, so $K$ is a functor and $\varkappa_c$ is a component of natural transformation $\varkappa:\ K\rightharpoonup T$. Then the diagram of direct sum~\eqref{dirsum-c} yields a diagram in $\mathcal{A}$
$$\xymatrix{
0\ar[r] & K(c_2) \ar[r]^{\varkappa_{c_2}} \ar@<1ex>[d]^-{K(\iota_2)}& T(c_2) \ar[r]^{\varphi_{c_2}} \ar@<1ex>[d]^-{T(\iota_2)}& S(c_2)\ar@<1ex>[d]^-{S(\iota_2)} \\
0\ar[r] & K(c_1) \oplus K(c_2) \ar[r]^{\varkappa_{c_1 \times c_2}} \ar@<-1ex>[d]_-{K(\pi_1)} \ar@<1ex>[u]^-{K(\pi_2)} & T(c_1) \oplus T(c_2) \ar[r]^{\varphi_{c_1 \times c_2}} \ar@<-1ex>[d]_-{T(\pi_1)} \ar@<1ex>[u]^-{T(\pi_2)} & S(c_1) \oplus S(c_2)\ar@<-1ex>[d]_-{S(\pi_1)} \ar@<1ex>[u]^-{S(\pi_2)} \\
0\ar[r] & K(c_1) \ar[r]^{\varkappa_{c_1}} \ar@<-1ex>[u]_-{K(\iota_1)}& T(c_1) \ar[r]^{\varphi_{c_1}} \ar@<-1ex>[u]_-{T(\iota_1)}& S(c_1)\ar@<-1ex>[u]_-{S(\iota_1)}
}$$
with exact rows and obvious commutativity properties. Since $\varkappa_c$ is a monomorphism, we see that left column defines a diagram of direct sum, so $K$ is again a product preserving functor.

Let $\lambda$ be a natural transformation such that $\lambda\varphi=0$, then for any object $c$ of $\mathcal{C}$ there exists a morphism ${\lambda'}_c$ in $\mathcal{A}$ such that $\lambda_c={\lambda'}_c \varkappa_c$ with $\varkappa$ constructed above. Since $\varkappa_c$ is a monomorphism, it can be easily proved that the above decomposition of $\lambda_c$ is functorial like the case of ordinary functors, so ${\lambda'}_c$ are components of a natural transformation. In other words, $\varkappa\in \ker \varphi$; dual statement may be proved analogously.  This verifies axioms (Abel-1) and (Abel-2) in $\mathop{\mathrm{FPFunct}} (\mathcal{C},\mathcal{A})$. Axiom (Abel-3) may be verified along the same lines as (Abel-1) was, so we have proved the following theorem.

\begin{theorem}
	Let $\mathcal{C}$ be a small additive category  and let $\mathcal{A}$ be an abelian category, then the category $\mathop{\mathrm{FPFunct}} (\mathcal{C},\mathcal{A})$ of product preserving functors from $\mathcal{C}$ to $\mathcal{A}$ is an abelian category with kernels and cokernels being defined objectwise.
\end{theorem}
  
Let $S$ be a set and let $S^\ast$ be ``the'' set of words generated by alphabet  $S$ (including an empty word), then  $S^\ast$ is a category with obvious finite products. Let $\mathbf{T}$ be an equationally defined algebraic theory of signature $S$ and let $\mathop{\mathrm{Alg}} (\mathbf{T},\mathcal{A})$ be the category of models of $\mathbf{T}$ in a category $\mathcal{A}$ with finite products. Then it is well known that there exist a small locally finite category $\mathcal{T}$ with finite products  such that  $\mathop{\mathrm{Alg}} (\mathbf{T},\mathcal{A})$ is equivalent to      $\mathcal{M}=\mathop{\mathrm{FPFunct}} (\mathcal{T},\mathcal {A})$ and we have a product preserving functor $T:\ S^\ast \rightarrow \mathcal {T}$ which is identity on objects~\cite{Ben,ATbook,Crole}.
 
Let $I$ be a discrete subcategory of $S^\ast$, then $\mathop{\mathrm{Funct}} (I,\mathcal {A})$ reduces to the category  $\mathcal{A}^I$  of functions from $I$ to  $\mathcal{A}$ with pointwise structure of abelian category inherited from $\mathcal{A}$. For any object $F$ of $\mathcal{M}$ define an object $
\Box_I F$ of $\mathcal{A}^I$ as a composition of $F$ with the inclusion $U_I: \  \xymatrix{I \ar@{^{(}->}[r] & S^\ast \ar[r]^T & \mathcal{T}}$, and natural transformations go to maps of functions. This defines an exact functor $\Box_I$ from $\mathcal{M}$ to $\mathcal{A}^I$, and we obtain the following proposition as a consequence of~\eqref{phi-times}, cf.~\cite{Ben}, Corollary 1.2.2.2. 

\begin{proposition}\label{prop-faith}
Suppose $S\subset I$, then $\Box_I$ is faithful.
\end{proposition}

Suppose now that $\mathcal{A}$ is complete, this is true e.g.  for $\mathcal{A}=\mathop{\mathrm{Alg}} (\mathbf{T},\mathcal{A}')$ provided $\mathcal{A}'$ is. Let $\boxtimes$ be tensor product of preadditive categories. Then it was shown in~\cite{Mitch}  that for a preadditive category $\mathcal{I}$ and a complete abelian category $\mathcal{A}$ there exists a preadditive bifunctor 
$$\mathrm{hom}_\mathcal{I} (-,-): \ (\mathop{\mathrm{Funct}} (\mathcal{I}, \mathbf{Ab}))^\mathrm{op} \boxtimes \mathop{\mathrm{Funct}} (\mathcal{I},\mathcal{A}) \rightarrow \mathcal{A}$$
  limit preserving in the first variable, and   it may be extended to preadditive functor $\mathop{\mathrm{Funct}} (\mathcal{C}, \mathop{\mathrm{Funct}} (\mathcal{I}^\mathrm{op}, \mathbf{Ab})) \boxtimes \mathop{\mathrm{Funct}} (\mathcal{I},\mathcal{A}) \rightarrow \mathop{\mathrm{Funct}} (\mathcal{C},\mathcal{A})$ also denoted by $\mathrm{hom}_\mathcal{I} (-,-)$ (this is called ``picking up the operators'' in~\cite{Mitch}, pp.14-15) in such a way that $\mathrm{hom}_\mathcal{I} (F,G)$ is additive provided the first argument is additive considered as a functor from $\mathcal{C}$ to $\mathop{\mathrm{Funct}} (\mathcal{I}^\mathrm{op},\mathcal{A})$.

Consider now $\mathcal{C}$ as a functor from $\mathcal{C}^\mathrm{op}\otimes \mathcal{C}$ to $\mathbf{Ab}$. This gives the functor of Yoneda embedding $Y_\mathcal{C}(-): \ \mathcal{C^\mathrm{op}}\rightarrow \mathop{\mathrm{Funct}} (\mathcal{C},\mathbf{Ab})$, $c\mapsto \mathrm{Hom}_\mathcal{C}(c,-)$.  Let $U_\mathcal{I}:\ \mathcal{I}\rightarrow  \mathcal{C}$ be a preadditive functor and let $U^*_\mathcal{I}$ be the corresponding neglecting functor from $\mathrm{Funct} (\mathcal{C}, \mathbf{Ab})$ to $\mathrm{Funct} (\mathcal{I}, \mathbf{Ab})$.   Then the right adjoint to the neclecting functor $\Box_\mathcal{I}$ from $\mathop{\mathrm{Funct}} (\mathcal{C},\mathcal{A})$ to $\mathop{\mathrm{Funct}} (\mathcal{I},\mathcal{A})$ reads~\cite{Mitch}
	
	\[ \mathrm{Ran}_{U_\mathcal{I}}(-)=\mathrm{hom}_\mathcal{I} ((Y_\mathcal{C})^\mathrm{op} U^*_\mathcal{I}, -) .\]
 
Dually, suppose that $\mathcal{A}$ is cocomplete, then there exist the bifunctor 	

$$ -\otimes_\mathcal{I} -:\   \mathop{\mathrm{Funct}} (\mathcal{I},  \mathcal{A}) \boxtimes \mathop{\mathrm{Funct}} (\mathcal{C}, \mathop{\mathrm{Funct}} (\mathcal{I}^\mathrm{op}, \mathbf{Ab})) \rightarrow \mathop{\mathrm{Funct}} (\mathcal{C},\mathcal{A})$$	
colimit preserving in the second variable  and the left adjoint to $\Box_\mathcal{I}$ has the form

\[ \mathrm{Lan}_{U_\mathcal{I}} (-)= - \otimes_\mathcal{I} Y_\mathcal{C}^\mathrm{op} (U^*_\mathcal{I})^\mathrm{op}, \]
where $Y_\mathcal{C}^\mathrm{op}$ is the product-preserving ``second Yoneda embedding'' $\mathcal{C}\rightarrow \mathop{\mathrm{Funct}} (\mathcal{C}^\mathrm{op},\mathbf{Ab})$ given by $c\mapsto \mathrm{Hom}_\mathcal{C}(-,c)$ (cf. \cite{Barr}, Prop.~5.3.18). Observe that  $Y_\mathcal{C}$ preserves zero object, then putting $\mathcal{I}=I$ and using the identification~\eqref{fp-zn} we obtain the following proposition.

	\begin{proposition}\label{prop-adj}
		Suppose  $\mathcal{A}$ is  cocomplete, then the functor $\Box_I:\ \mathcal{M}\rightarrow \mathcal{A}^I$ has  left adjoint. 
	\end{proposition}

	\section{Cohomology}
	
 Combining Prop.~\ref{prop-faith} and Prop.~\ref{prop-adj}, we can construct   a resolvent pair of abelian categories in the sense of~(\cite{H}, Ch. IX \S\S 5,6) for any $I \supset S$ using the pair of adjoint functors $L_{\Box_I}\dashv  \Box_I$, so we can define relative extension functor $\mathop{\mathrm{Ext}}^\ast_{\Box_I} (-,-)$ from $\mathcal{M}^{\mathrm{op}}\times \mathcal {M}$ to $\mathbf{Ab}$ w.r.t. the proper  class of $\Box_I$-split short exact sequences along the lines of (\cite{H}, Ch. XII \S\S 4,5), and we obtain the following proposition.

\begin{proposition}\label{prop-ext}
	Suppose  $\mathcal{A}$ is complete, then an isomorphism
		\[ H^n_{\Box_I} (F,G):=\mathrm{Ext}^n_{\Box_I} (F,G) \cong H^n (\mathrm{Hom}_\mathcal{M} (\beta^{\Box_I} (F), G) \]
		there exist, where $\beta^{\Box_I} (F)$  is a relatively projective resolvent of $F$ in $\mathcal{M}$ constructed using the pair of adjoint functors $\mathrm{Lan}_{\Box_I}\dashv \Box_I $.
	\end{proposition}

Then proceeding along the lines of~\cite{CE}, Ch. XVI \S 1 and \cite{H}, Ch.XII \S\S9,10 we obtain the following proposition which relates this  cohomology to "absolute" one, thus extending  Theorem~C of~\cite{JP91}.

\begin{proposition}
Suppose  $\mathcal{A}$ is both complete and cocomplete, then an isomorphism
\[ H^n_{\Box_I} (F,G) \cong H^n(\mathrm{Hom}_{\mathop{\mathrm{Funct}}(\mathcal{T},\mathcal{A})} ( F,\beta_{\Box_I} (G)) \eqno (*)\]
there exist, where $\beta_{\Box_I}$ is a relatively injective coresolvent of $G$ in $\mathop{\mathrm{Funct}}(\mathcal{T},\mathcal{A})$ constructed using the pair of adjoint functors $\Box_I\dashv \mathrm{Ran}_{\Box_I} $. Since $\mathop{\mathrm{Funct}}(\mathcal{T},\mathcal{A})$ has enough "absolute" injectives, then we can take "absolute" injective resolvent of $G$ in (*) instead of relative one.
\end{proposition}

\begin{example}
Let $\mathcal{A}=\mathop{\mathrm{Funct}} (\mathcal{C}^\mathrm{op},\mathbf{Ab})$ and consider resolvent of the functor $Y_\mathcal{C}^\mathrm{op}$. Then like the case of rings and modules (\cite{H}, Ch.X, \S 2) we obtain that $\beta^{\Box_I}_{n+1} (Y_\mathcal{C}^\mathrm{op})$ is the functor from $\mathcal{C}$ to $\mathcal{A}$ defined by
$$c \mapsto \bigoplus_{i_1, \ldots, i_n \in I} \mathrm{Hom}_\mathcal{C}(-,i_1) \otimes_\mathbf{Z} \mathrm{Hom}_\mathcal{C}(i_1,i_2)  \otimes_\mathbf{Z} \ldots  \otimes_\mathbf{Z}  \mathrm{Hom}_\mathcal{C}(i_n,c),$$
where explicit formula for tensor product of functors (\cite{Mitch},p.26) was used. This resolvent coincides with the ordinary Hochschild-Mitchell one (\cite{Mitch}, p.70) except for the range of indices $i_1,\ldots i_n$ which run  $I$ but not over the whole $\mathop{\mathrm{Ob}} \mathcal{C}$.

Let $F$ be an element of $\mathcal{M}=\mathop{\mathrm{FPFunct}} (\mathcal{C},\mathcal{A})$ and let $\phi (f,f_1,\ldots,f_n,f')$ be $\mathbf{Z}$-linear function from 
$$\prod_{i,i_1,\ldots,i_n,i'} \mathrm{Hom}_\mathcal{C} (i,i_1) \times \mathrm{Hom}_\mathcal{C} (i_1,i_2) \times \ldots \times \mathrm{Hom}_\mathcal{C} (i_n,i')$$ 
to $\mathrm{Hom}_\mathbf{Ab}(\mathrm{Hom}_\mathcal{C} (i,i'), F(i,i'))$
 with differential  
\begin{multline}
\delta\phi (f,f_1,\ldots,f_n,f')=F(1\boxtimes f) \phi (f_1,\ldots,f_n,f') +\\ 
+\sum_{p=1}^n (-1)^p \phi(f,f_1,\ldots, f_pf_{p+1}, \ldots,f_n,f') + \\
+ \phi (f,f_1,\ldots,f_n)F(f'\boxtimes 1),
\end{multline}
where $F$ is considered as a functor from $\mathcal{C}\boxtimes \mathcal{C}^\mathrm{op}$ to $\mathbf{Ab}$. Denote of this subcomplex with $I=S$ by $H^n_\mathrm{red}(Y_\mathcal{C}^\mathrm{op},F)$, then using~(\ref{phi-times}) we obtain the following decomposition of full cohomology groups
$$H^n(Y_\mathcal{C}^\mathrm{op},F)=\bigoplus_{S^\ast \times S^\ast} H^n_\mathrm{red}(Y_\mathcal{C}^\mathrm{op},F).$$

\end{example}

\end{document}